\documentclass[12pt,a4paper]{article}

\usepackage{amsmath,amssymb}

\newcommand{\rset}{\mathbf{R}}
\newcommand{\rmax}{\rset_{\max}}
\newcommand{\smaxmin}{S_{\max,\min}}
\newcommand{\Mat}{\mathrm{Mat}}
\newcommand{\0}{\mathbf{0}}
\newcommand{\1}{\mathbf{1}}
\newcommand{\cle}{\preccurlyeq}

\newcommand{\aint}{\mathbf a}
\newcommand{\A}{\mathbf A}
\newcommand{\B}{\mathbf B}
\newcommand{\x}{\mathbf x}
\newcommand{\X}{\mathbf X}
\newcommand{\y}{\mathbf y}
\newcommand{\lA}{\underline{\A}}
\newcommand{\la}{\underline{\aint}}
\newcommand{\lB}{\underline{\B}}
\newcommand{\lx}{\underline{\x}}
\newcommand{\ly}{\underline{\y}}
\newcommand{\uA}{\overline{\A}}
\newcommand{\ua}{\overline{\aint}}
\newcommand{\uB}{\overline{\B}}
\newcommand{\ux}{\overline{\x}}
\newcommand{\uy}{\overline{\y}}

\newtheorem{thm}{Theorem}
\newtheorem{predl}{Proposition}

\begin{document}

\small

Doklady Mathematics, vol.~62, No.~2, 2000, pp. 199--201

Translated from Doklady Akademii Nauk, vol. 374, No. 3, 2000, pp. 304--306

\bigskip

\large

\begin{center}
G.~L.~Litvinov, A.~N.~Sobolevskii\\
\bigskip
EXACT INTERVAL SOLUTIONS TO THE DISCRETE BELLMAN EQUATION\\
AND POLYNOMIAL COMPLEXITY OF PROBLEMS\\
IN INTERVAL IDEMPOTENT LINEAR ALGEBRA
\end{center}

In this paper, we construct a solution to a linear matrix interval equation
of the form $X = \A X + \B$ (the discrete stationary Bellman equation) over
partially ordered semirings, including the semiring $\rset_+$ of
nonnegative real numbers and all idempotent semirings. We also discuss the
computational complexity of problems in interval idempotent linear algebra
(for more detail on idempotent mathematics, see, e.g., \cite{mk,bsy}). In
traditional interval analysis, problems of this kind are generally
$NP$-hard \cite{klrk,sh}. In this paper, we consider matrix equations over
positive semirings (in the sense of \cite{golan}); in this case the
computational complexity of the problem is polynomial.

Idempotent and other positive semirings naturally arise in optimization
problems. Many of these problems turn out to be linear over appropriate
idempotent semirings \cite{mk,bsy}. In this case, the system of equations
$X = \A X + \B$ is a natural analogue of a usual linear system in
traditional linear algebra over fields. Carr\'e~\cite{carre} showed that
many of the well-known algorithms of discrete optimization are analogous to
standard algorithms in traditional computational linear algebra.

1. Consider a semiring, i.e., a set~$S$ endowed with two associative
operations, addition~$\oplus$ and multiplication~$\odot$, such that
addition is commutative, multiplication is distributive over addition from
either side, $\0$ and $\1$ are the respective neutral elements of addition
and multiplication, $\0 \odot x = x \odot \0 = \0$ for all $x \in S$, and
$\0 \neq \1$. Let the semiring~$S$ be partially ordered by a
relation~$\cle$ such that $\0$ is the least element and the inequality $x
\cle y$ implies that $x \oplus z \cle y \oplus z$, $x \odot z \cle y \odot
z$, and~$z \odot x \cle z \odot y$ for all $x, y, z \in S$; in this case
the semiring~$S$ is called positive (see, e.g., \cite{golan}).

A semitring~$S$ is called idempotent if $x \oplus x = x$ for all $x \in S$
\cite{mk,bsy,danlm}. Addition~$\oplus$ defines a canonical partial
order~$\stackrel\oplus\cle$ on~$S$ by the rule $x\stackrel\oplus\cle y$ iff
$x \oplus y = y$.  Any idempotent semiring is positive with respect to this
order.  Note also that $x \oplus y = \sup\{x,y\}$ with respect to the
canonical order.  In what follows, we assume that all idempotent semirings
are ordered by the canonical partial order relation.

The best known and most important examples of positive semirings are
``numerical'' semirings consisting of (a subset of) real numbers and
ordered by the conventional linear order $\leqslant$ on~$\rset$: the
semiring~$\rset_+$ with the usual operations $\oplus = +$, $\odot = \cdot$
and neutral elements $\0 = 0$, $\1 = 1$; the semiring~$\rmax = \rset \cup
\{-\infty\}$ with the operations $\oplus = \max$, $\odot = +$ and neutral
elements $\0 = -\infty$, $\1 = 0$; the semiring $\hat\rmax = \rmax \cup
\{\infty\}$, where $x \cle \infty$, $x \oplus \infty = \infty$ for all $x$,
$x \odot \infty = \infty \odot x = \infty$ if $x \neq \0$, and $\0 \odot
\infty = \infty \odot \0$; and the semiring~$\smaxmin^{[a,b]} = [a, b]$,
where $-\infty \leqslant a < b \leqslant +\infty$, with the operations
$\oplus = \max$, $\odot = \min$ and neutral elements $\0 = a$, $\1 = b$.
The semirings~$\rmax$, $\hat\rmax$, and~$\smaxmin^{[a,b]} = [a, b]$ are
idempotent.

Many mathematical constructions, concepts, and results over the fields of
real and complex numbers have nontrivial analogues over idempotent
semirings.  Idempotent semirings have recently become the subject of a new
branch of mathematics, idempotent analysis \cite{mk,bsy,danlm}.

Let a positive semiring~$S$ be endowed with a partial unary closure
operation~$*$ such that $x \cle y$ implies $x^* \cle y^*$ and $x^* = \1
\oplus (x^* \odot x) = \1 \oplus (x \odot x^*)$ on its domain of
definition. In particular, $\0^* = \1$ by definition. These axioms imply
that $x^* = \1 \oplus x \oplus x^2 \oplus \dots \oplus (x^* \odot x^n)$ if
$n \geqslant 1$. Thus $x^*$ can be considered as a `regularized sum' of the
series $x^* = \1 \oplus x \oplus x^2 \oplus \dots$; in an idempotent
semiring, by definition, $x^* = \sup \{ \1, x, x^2, \dots \}$ if this
supremum exists.

In numerical semirings, the operation~$*$ is defined as follows: $x^* = (1
- x)^{-1}$ if $x \prec 1$ in $\rset_+$, $x^* = \1$ if $x \cle \1$ in
$\rmax$ and $\hat\rmax$, $x^* = \infty$ if $x \succ \1$ in $\hat\rmax$,
$x^* = \1$ for all $x$ in $\smaxmin^{[a,b]}$. In all other cases, $x^*$ is
undefined.  Note that the operation~$*$ is defined everywhere in idempotent
semirings that are $a$-complete in the sense of \cite{danlm} (e.g., in
$\hat\rmax$ or $\smaxmin^{[a,b]}$). For more detail, see \cite{golan}.

2. Let $S$ be a set partially ordered by a relation $\cle$.  A closed
interval in~$S$ is a subset of the form $\x = [\lx, \ux] = \{\, x \in S
\mid \lx \cle x \cle \ux\, \}$, where the elements $\lx \cle \ux$ are
called lower and upper bounds on the interval $\x$.  The order~$\cle$
induces a partial ordering on the set of all closed intervals in~$S$: $\x
\cle \y$ iff $\lx \cle \ly$ and $\ux \cle \uy$.

A weak interval extension $I(S)$ of a positive semiring~$S$ is the
set of all closed intervals in~$S$ endowed with the operations $\oplus$
and~$\odot$ defined by ${\x \oplus \y} = [{\lx \oplus \ly}, {\ux \oplus
\uy}]$, ${\x \odot \y} = [{\lx \odot \ly}, {\ux \odot \uy}]$ and with a
partial order induced by the order in $S$. The closure operation in $I(S)$
is defined by $\x^* = [\lx^*, \ux^*]$ (see also \cite{golan}; for interval
analysis over $\rset$, see, e.g., \cite{moore}).

\begin{predl}
The weak interval extension $I(S)$ of a positive semiring~$S$ is closed
under the operations $\oplus$ and~$\odot$ and forms a positive semiring
with a zero element $[\0, \0]$ and a unit element $[\1, \1]$. The interval
${\x \oplus \y}$ $({\x \odot \y})$ contains the set $\{\, {x \oplus  y} \mid
x \in \x, y \in \y \,\}$ $(\{\, {x \odot  y} \mid x \in \x, y \in \y \,\}$,
respectively$)$ and its bounds are elements of this set.
\end{predl}

3. Denote by $\Mat_{mn}(S)$ a set of all matrices $A = (a_{ij})$ with
$m$~rows and $n$~columns, whose coefficients belong to a semiring~$S$.  The
sum $A \oplus B$ of matrices $A, B \in \Mat_{mn}(S)$ and the product $AB$
of matrices $A \in \Mat_{lm}(S)$ and $B \in \Mat_{mn}(S)$ are defined
according to the usual rules of linear algebra. If the semiring~$S$ is
positive, then the set $\Mat_{mn}(S)$ is ordered by the relation $A =
(a_{ij}) \cle B = (b_{ij})$ iff $a_{ij} \cle b_{ij}$ in~$S$ for all $1
\leqslant i \leqslant m$, $1 \leqslant j \leqslant n$.

Matrix multiplication is consistent with the order~$\cle$ in the following
sense: if $A, A' \in \Mat_{lm}(S)$, $B, B' \in \Mat_{mn}(S)$ and $A \cle
A'$, $B \cle B'$, then $AB \cle A'B'$ in $\Mat_{ln}(S)$. The set
$\Mat_{nn}(S)$ of square matrices of order~$n$ over a (positive,
idempotent) semiring~$S$ forms a (positive, idempotent) semiring with a
zero element $O = (o_{ij})$, where $o_{ij} = \0$, $1 \leqslant i, j
\leqslant n$, and with a unit element $E = (\delta_{ij})$, where
$\delta_{ij} = \1$ if $i = j$ and $\delta_{ij} = \0$ in the opposite case.

The closure operation in matrix semirings over a positive semiring~$S$ can
be  defined inductively (for another way of doing this, see~\cite{golan}):
$A^* = (a_{11})^* = (a^*_{11})$ in $\Mat_{11}(S)$ and for any integer $n >
1$ and any matrix
$$
   A = \begin{pmatrix} A_{11}& A_{12}\\ A_{21}& A_{22} \end{pmatrix},
$$
where $A_{11} \in \Mat_{kk}(S)$, $A_{12} \in \Mat_{k\, n - k}(S)$,
$A_{21} \in \Mat_{n - k\, k}(S)$, $A_{22} \in \Mat_{n - k\, n - k}(S)$,
$1 \leqslant k \leqslant n$, by definition,
$$
   A^* = \begin{pmatrix}
   A^*_{11} \oplus A^*_{11} A_{12} D^* A_{21} A^*_{11} &
   A^*_{11} A_{12} D^* \\
   D^* A_{21} A^*_{11} &
   D^*
   \end{pmatrix},
$$
where $D = A_{22} \oplus A_{21} A^*_{11} A_{12}$. It can be proved that
this definition of $A^*$ implies that $A^* = A^*A \oplus E$, and, thus,
$A^*$ is a ``regularized sum'' of the series $E \oplus A \oplus A^2 \oplus
\dots$.

Note that this recurrence relation coincides with the formulas of the
escalator method for matrix invertion in traditional linear algebra over
the field of real or complex numbers up to the algebraic operations used.
Hence this algorithm of matrix closure is polynomial in~$n$.

Let $S$ be a positive semiring and $\A = (\aint_{ij}) \in \Mat_{mn}(I(S))$
be a matrix whose coefficients are closed intervals in~$S$. The matrices
$L(\A) = (\la_{ij}), U(\A) = (\ua_{ij}) \in \Mat_{mn}(S)$ are called the
lower and the upper matrices of the interval matrix $\A$.  Evidently,
$L(\A) \cle U(\A)$ in $\Mat_{mn}(S)$.

Since, for any positive semiring~$S$, the sets $I(S)$ and $\Mat_{nn}(S)$
form positive semirings, the sets $I(\Mat_{nn}(S))$ and~$\Mat_{nn}(I(S))$
form positive semirings with respect to the operations defined above.
\begin{predl}
The semirings $I(\Mat_{nn}(S))$ and~$\Mat_{nn}(I(S))$ are isomorphic to
each other, and the isomorphism is defined by $\A \in \Mat_{nn}(I(S))
\mapsto [L(\A), U(\A)] \in I(\Mat_{nn}(S))$.
\end{predl}

By definition, the addition and multiplication of matrix intervals
in~$I(\Mat_{nn}(S))$ are reduced to separate matrix addition and
multiplication of their lower and upper matrices. An analogous statement
for lower and upper matrices in $\Mat_{nn}(I(S))$ follows from the last
proposition.

4. Let~$S$ be a positive semiring. The discrete stationary Bellman equation
has the form
\begin{equation*}
	X = AX \oplus B, \tag*{$(*)$}
\end{equation*}
where $A \in \Mat_{nn}(S)$, $X, B \in \Mat_{ns}(S)$, and the matrix~$X$ is
unknown. Let $A^*$ be the closure of the matrix~$A$. It follows from the
identity $A^* = A^*A \oplus E$ that the matrix $A^*B$ satisfies this
equation; moreover, it can be proved that, in idempotent semirings, this
solution is the least in the set of solutions to equation~$(*)$ with
respect to the partial order in $\Mat_{ns}(S)$.

Let $A = \A \in \Mat_{nn}(I(S))$, $B = \B \in \Mat_{ns}(I(S))$. The unified
(least) solution set of equation~$(*)$ is the set $\Sigma(\A,\B) = \{\,
A^*B \mid A \in \A, B \in \B \,\}$. The interval $\X = \A^*\B \in
\Mat_{ns}(I(S))$ that satisfies equation~$(*)$ in the algebraic sense is
called the (least) algebraic solution to this equation. Other defintions of
a solution set of a matrix interval linear equation can be found, for
example, in  \cite{sh}.

Our main result is the following theorem.
\begin{thm}
The closed interval $[L(\A^*\B), U(\A^*\B)]$ in $\Mat_{ns}(S)$ that
corresponds to an algebraic solution to equation~$(*)$ contains the unified
solution set $\Sigma(\A,\B)$ of equation~$(*)$, and the bounds of this
interval belong to $\Sigma(\A,\B)$. The algebraic solution $\A^*\B$ can be
constructed in a polynominal number of operations.
\end{thm}

The proof follows from the fact that matrix multiplication and closure are
consistent with the partial order in a matrix semiring. The lower and the
upper matrices of an algebraic solution $\A^*\B$ to equation~$(*)$ satisfy
the point equations $X = \lA X \oplus \lB$ and $X = \uA X \oplus \uB$, and
algebraic solutions to these equations can be constructed by the matrix
closure algorithm described in section~3, which is polynomial in~$n$.

Note that this theorem was proved in the paper \cite{bn} (see also
\cite{ah}, Theorem~12.2) in the case of interval linear algebra over the
semiring~$\rset_+$ of nonnegative real numbers.

Under some natural additional conditions on the operations $\oplus$,
$\odot$, and $*$
a stronger equality $[L(\A^*\B), U(\A^*\B)] = \Sigma(\A,\B)$ holds
in the case of an idempotent semiring~$S$.
Moreover, as far as we know, there are no $NP$-hard
computational problems in interval linear algebra in this case. This is
consistent with the general observation that idempotent analogues of
constructions in traditional mathematics over numerical fields are
considerably simpler than their prototypes \cite{mk,bsy,danlm}.

ACKNOWLEDGEMENTS: The authors are grateful to V.~P.~Maslov and S.~P.~Shary
for useful discussions.

The work was supported by Russian Foundation for Basic
Research (project no.~99--01--01198).


\bigskip

E--mail: litvinov@islc.msk.su

glitvinov@mail.ru

\end{document}